\documentclass[12pt]{article}
\usepackage{amsmath}
\usepackage{amsfonts}
\usepackage{amsthm}

\newtheorem{lemma}{Lemma}[section]
\newtheorem{theorem}[lemma]{Theorem}
\newtheorem{corollary}[lemma]{Corollary}

\title{Rotations in three, four, and five dimensions}
\author{jason hanson}

\begin{document}
\maketitle

\begin{abstract}
The geometry of rotations in dimensions 3, 4, and 5 is discussed using the matrix exponential map.  Explicit closed formulas for the exponential of an antisymmetric matrix, as well as the logarithm of a rotation, are given for these dimensions.
\end{abstract}

\section{Introductory remarks}

For representing rotations in three dimensions, one usually considers two choices: Euler angles, or quaternions.  For aerospace navigation, Euler angles are a natural choice, since the familiar maneuvers of yaw, pitch, and roll are coordinate axis rotations.  However, from a global point of view (as opposed to the local viewpoint of a pilot), Euler angles are unnatural and awkward.  In particular, the decomposition of a rotation about a noncoordinate axis into Euler angles is not unique and depends on the sequence of axis rotations.  On the other hand, quaternions can be used to represent a rotation in a more convenient manner: rotation about an arbitrary axis can be written in a simple form, and given a rotation, the extraction of the rotation axis and angle is a relatively straight--forward process.  The price paid for this convenience, of course, is the introduction of a less familiar Mathematical apparatus.

In less common use is a third choice: the representation of a rotation by a linear map, or matrix.  Fundamentaly, the Euler angle representation is a matrix representation, but one only uses the representation of coordinate axis rotations, which have a very simple form.  In some texts, one may encounter a matrix formula for rotation about an arbitrary axis, but usually this is presented in the form of a rather unappealing and unyieldy $3\times 3$ matrix.  In this article, we will recast this matrix in a more appealing form by using some simpler linear maps.  As a natural byproduct of this approach, we will obtain a simple way of extracting the rotation axis and angle, but using a more familiar apparatus than quaternions.  Moreover, we will be able to extend the constructions to higher dimensions.  While one may use quaternions to represent rotations in four dimensions (see \cite{Weiner}), one cannot go beyond.

\section{Rotations in three dimensions}

A rotation $R$ in three dimensional space is specified by a (real) scalar, the {\em rotation angle} $\theta$, and a (real) unit vector ${\bf u}$, the {\em rotation axis}.  The effect of $R$ on any vector ${\bf v}$ is given by {\em Rodrigues' Rotation Formula:}
\begin{equation}\label{eq:rodrigues}
  R{\bf v}
  =\cos\theta\,{\bf v}
    +(1-\cos\theta)({\bf u}\cdot{\bf v}){\bf u}
    +\sin\theta\,{\bf u}\times{\bf v}.
\end{equation}
Geometrically, if we write ${\bf v}$ in terms of its components parallel and perpendicular to ${\bf u}$ (respectively): ${\bf v}={\bf v}_\parallel+{\bf v}_\perp$, then $R{\bf v}_\parallel={\bf v}_\parallel$, and $R{\bf v}_\perp$ is the rotation of ${\bf v}_\perp$ in the plane orthogonal to ${\bf u}$ by $\theta$ radians.  It is important to note that this is {\em counterclockwise} rotation about ${\bf u}$; indeed, using the standard three--dimensional vector identity ${\bf u}\times({\bf v}\times{\bf w})=({\bf u}\cdot{\bf w}){\bf v}-({\bf u}\cdot{\bf v}){\bf w}$, one computes ${\bf v}_\perp\times({\bf u}\times{\bf v})=|{\bf v}_\perp|^2{\bf u}$, which is a positive multiple of ${\bf u}$.

\subsection{Rotation as a linear map}

We may express equation \eqref{eq:rodrigues} as a $3\times 3$ matrix, or equivalently, a linear map.  Let ${\it proj}_{\bf u}$ denote the linear map ${\it proj}_{\bf u}({\bf v})\doteq({\bf u}\cdot{\bf v}){\bf u}$ (we are still assuming that $|{\bf u}|=1$); that is, ${\it proj}_{\bf u}$ is orthogonal projection onto ${\bf u}$.  In addition, let $\Lambda_{\bf u}$ denote the linear map $\Lambda_{\bf u}({\bf v})\doteq{\bf u}\times{\bf v}$.  Explicitly as matrices,
\begin{equation*}
  {\it proj}_{\bf u}
  =\begin{pmatrix} u_1^2  & u_1u_2 & u_1u_3\\
                   u_1u_2 & u_2^2  & u_2u_3\\
                   u_1u_3 & u_2u_3 & u_3^2
   \end{pmatrix}
  \quad\text{and}\quad
  \Lambda_{\bf u}
  =\begin{pmatrix}   0 & -u_3 &  u_2\\
                   u_3 &    0 & -u_1\\
                  -u_2 &  u_1 &    0
   \end{pmatrix}
\end{equation*}
where ${\bf u}=(u_1,u_2,u_3)$.  Using these maps, equation \eqref{eq:rodrigues} becomes
\begin{equation}\label{eq:rot3}
  R=\cos\theta\,I
    +(1-\cos\theta)\,{\it proj}_{\bf u}
    +\sin\theta\,\Lambda_{\bf u}.
\end{equation}
Here $I$ denotes the identity map: $I{\bf v}\doteq{\bf v}$.

\subsection{Finding the angle and axis of a rotation}

Conversely, given a rotation $R$, we may deduce its rotation angle and axis.  For the angle, we take the trace of $R$ in equation \eqref{eq:rot3} to get ${\it tr}(R)=1+2\cos\theta$.  Thus
\begin{equation}\label{eq:angle3}
  \theta
  =\cos^{-1}\tfrac{1}{2}\bigl({\it tr}(R)-1\bigr).
\end{equation}
For the axis, we make use of the fact that $I$ and ${\it proj}_{\bf u}$ are symmetric, while $\Lambda_{\bf u}$ is antisymmetric, so that equation \eqref{eq:rot3} also yields
\begin{equation}\label{eq:axis3}
  \tfrac{1}{2}(R-R^t)=\sin\theta\,\Lambda_{\bf u},
\end{equation}
where $R^t$ denotes the tranpose of $R$.  The components of ${\bf u}$ can be read directly off of $\Lambda_{\bf u}=(R-R^t)/2\sin\theta$ (or more practically, since $2\sin\theta\,\Lambda_{\bf u}=\Lambda_{2\sin\theta\,{\bf u}}$, we can just read off the components of $2\sin\theta\,{\bf u}$ from $R-R^t$, and then normalize to get ${\bf u}$).

Two remarks are in order.  First, the rotation angle $\theta$ determined by equation \eqref{eq:angle3} is in the range $0\leq\theta\leq\pi$, whereas there is no such restriction in equation \eqref{eq:rot3}.  If the original rotation angle is greater than $\pi$ (but less that $2\pi$), then the axis determined from \eqref{eq:axis3} will be opposite of the original rotation axis.  Second, we can only use equation \eqref{eq:axis3} to determine the rotation axis if $0<\theta<\pi$.  In the case $\theta=0$, the rotation is the identity.  And in the case $\theta=\pi$, the rotation is reflection in the plane orthogonal to the vector ${\bf u}$, which we can deduce up to sign using ${\it proj}_{\bf u}=\tfrac{1}{2}(I+R)$.

\subsection{Rotation as an exponential}
\label{subsec:rotexp}

Recall that the Taylor series expansion of the exponential function is given by $e^x=\sum_{k\geq 0}x^k/k!$, which converges for all real numbers $x$.  Similary, the exponential of a square matrix $A$ is {\em defined} by the formula
\begin{equation*}
  \exp(A)
  \doteq\sum_{k=0}^\infty\frac{1}{k!}\,A^k.
\end{equation*}
Here we set $A^0\doteq I$.  It can be shown that this summation converges for any matrix square $A$ (see for example, \cite{Curtis}).

It is known from the theory of Lie groups that every rotation is the exponential of an antisymmetric matrix.  We can verify this directly in three dimensions.  Let ${\bf u}$ be a three--dimensional unit vector.  Starting from the matrix representation of $\Lambda_{\bf u}$, or from vector identities, one computes that $\Lambda_{\bf u}^2=-(I-{\it proj}_{\bf u})$ and $\Lambda_{\bf u}^3=-\Lambda_{\bf u}$.  In particular, $-\Lambda_{\bf u}^2$ is orthogonal projection onto the plane orthogonal to ${\bf u}$.  It follows that $\Lambda_{\bf u}^{2p}=(-1)^p(I-{\it proj}_{\bf u})$ for all $p\geq 1$, and $\Lambda_{\bf u}^{2p+1}=(-1)^p\Lambda_{\bf u}$ for all $p\geq 0$.  We then compute
\begin{align*}
  \exp(\theta\Lambda_{\bf u})
  &=\sum_{k=0}^\infty\frac{\theta^k}{k!}\,\Lambda_{\bf u}^k
   =I+\sum_{p=1}^\infty\frac{\theta^{2p}}{(2p)!}\,
                    \Lambda_{\bf u}^{2p}
    +\sum_{p=0}^\infty\frac{\theta^{2p+1}}{(2p+1)!}\,
                   \Lambda_{\bf u}^{2p+1}\\
  &=I+\left[\sum_{p=1}^\infty\frac{(-1)^p}{(2p)!}\,\theta^{2p}\right]
                          (I-{\it proj}_{\bf u})
    +\left[\sum_{p=0}^\infty\frac{(-1)^p}{(2p+1)!}\,\theta^{2p+1}\right]
                         \Lambda_{\bf u}\\
  &=I+(\cos\theta-1)(I-{\it proj}_{\bf u})+\sin\theta\,\Lambda_{\bf u}
\end{align*}
by virtue of the Taylor series expansions for sine and cosine.  This last expression is just a rearrangement of equation \eqref{eq:rot3}.  That is, $\exp(\theta\Lambda_{\bf u})$ is counterclockwise rotation by $\theta$ radians about ${\bf u}$.

\section{Simple rotations}

We now extend the above constructions used in three dimensions to any dimension.  To do this, we define two special linear maps.  If ${\bf u}$ and ${\bf v}$ are two vectors in $n$--dimensional space, their {\bf outer product} is the linear map ${\bf u}\otimes{\bf v}$ such that
\begin{equation*}
  ({\bf u}\otimes{\bf v})({\bf w})
  =({\bf v}\cdot{\bf w}){\bf v}
\end{equation*}
for all vectors ${\bf w}$.  This is a generalization of the projection operation: if $|{\bf u}|=1$, then ${\it proj}_{\bf u}={\bf u}\otimes{\bf u}$.  The {\bf wedge product} of ${\bf u}$ and ${\bf v}$ is defined as
\begin{equation*}
  {\bf u}\wedge{\bf v}
  \doteq{\bf u}\otimes{\bf v}-{\bf v}\otimes{\bf u}.
\end{equation*}
The wedge product will serve as our generalization of $\Lambda_{\bf u}$; in fact in three dimensions, we have ${\bf u}\wedge{\bf v}=-\Lambda_{{\bf u}\times{\bf v}}$.

\subsection{Fundamental computations}

The following properties of the outer and wedge products are readily verified from the definitions.

\begin{lemma}\label{lem:elemcomp}
For all vectors ${\bf a}$, ${\bf b}$, ${\bf u}$, ${\bf v}$, we have
\begin{enumerate}
\item ${\it tr}({\bf u}\otimes{\bf v})={\bf u}\cdot{\bf v}$
\item\label{it:outercomp} $({\bf a}\otimes{\bf b})({\bf u}\otimes{\bf v})=({\bf b}\cdot{\bf u})({\bf a}\otimes{\bf v})$.
\item $|{\bf u}\wedge{\bf v}|^2{\it proj}_{\bf uv}=|{\bf v}|^2{\bf u}\otimes{\bf u}-({\bf u}\cdot{\bf v})({\bf u}\otimes{\bf v}+{\bf v}\otimes{\bf u})+|{\bf u}|^2{\bf v}\otimes{\bf v}$.
\item\label{it:proj} $({\bf u}\wedge{\bf v})^2=-|{\bf u}\wedge{\bf v}|^2{\it proj}_{\bf uv}$
\item\label{it:cube} $({\bf u}\wedge{\bf v})^3=-|{\bf u}\wedge{\bf v}|^2{\bf u}\wedge{\bf v}$
\end{enumerate}
where $|{\bf u}\wedge{\bf v}|\doteq\sqrt{|{\bf u}|^2|{\bf v}|^2-({\bf u}\cdot{\bf v})^2}$, and ${\it proj}_{\bf uv}$ is orthogonal projection onto the plane spanned by ${\bf u}$ and ${\bf v}$.\qed
\end{lemma}

We remark that, as the notation indicates, $|{\bf u}\wedge{\bf v}|$ is the matrix norm of ${\bf u}\wedge{\bf v}$.  For antisymmetric (or symmetric) matrices $A$ and $B$, one usually defines their inner product as $A\cdot B\doteq\tfrac{1}{2}{\it tr}(AB^t)$.  In this sense, $|{\bf u}\wedge{\bf v}|^2=({\bf u}\wedge{\bf v})\cdot({\bf u}\wedge{\bf v})$.  In particular, ${\bf u}\wedge{\bf v}=0$ if and only if $|{\bf u}\wedge{\bf v}|=0$.  Note that we may also write $|{\bf u}\wedge{\bf v}|=|{\bf u}|\,|{\bf v}|\sin\phi$, where $\phi$ is the angle between ${\bf u}$ and ${\bf v}$.

Using lemma \ref{lem:elemcomp}\eqref{it:proj} and \eqref{it:cube}, we may again sum the Taylor series as we did in section \ref{subsec:rotexp} to compute the exponential of ${\bf u}\wedge{\bf v}$.

\begin{theorem}\label{thm:expsimple}
For all vectors ${\bf u}$, ${\bf v}$ with $|{\bf u}\wedge{\bf v}|=1$, and any real scalar $\theta$,
$$\exp(\theta\,{\bf u}\wedge{\bf v})
  =I-(1-\cos\theta)\,{\it proj}_{\bf uv}
    +\sin\theta\,{\bf u}\wedge{\bf v}.\qed
$$
\end{theorem}

\subsection{Rotation in a two--plane}

In any dimension greater than one, a (linear) {\em two--plane} is spanned by two linearly independent vectors.  We will call a linear map corresponding to rotation in a two--plane, and that fixes the orthogonal complement, a {\bf simple rotation}.

\begin{lemma}\label{lem:orient}
Given vectors ${\bf u}$, ${\bf v}$, if ${\bf a}=\alpha{\bf u}+\beta{\bf v}$ and ${\bf b}=\gamma{\bf u}+\delta{\bf v}$ for scalars $\alpha$, $\beta$, $\gamma$, $\delta$, then ${\bf a}\wedge{\bf b}=(\alpha\delta-\beta\gamma)\,{\bf u}\wedge{\bf v}$.
\end{lemma}

\begin{proof}
Use the linearity and antisymmetry of the wedge product.
\end{proof}

\begin{lemma}\label{lem:plane}
The vectors ${\bf u}$, ${\bf v}$ span a two--plane if and only if ${\bf u}\wedge{\bf v}\neq 0$.  Moreover, the vectors ${\bf a}$, ${\bf b}$ span the same two--plane if and only if ${\bf a}\wedge{\bf b}$ is a nonzero multiple of ${\bf u}\wedge{\bf v}$.
\end{lemma}

\begin{proof}
The first statement follows from the remarks after lemma \ref{lem:elemcomp}.  The necessity of the second follows from the previous lemma;  for the sufficiency, use lemma \ref{lem:elemcomp}\eqref{it:proj} to show that ${\it proj}_{\bf ab}={\it proj}_{\bf uv}$.
\end{proof}

An {\em orientation} for the two--plane spanned by ${\bf u}$, ${\bf v}$ is a choice of ordering of the vectors, say $({\bf u},{\bf v})$.  By lemma \ref{lem:orient}, a pair of vectors $({\bf a},{\bf b})$ from this two--plane has the same orientation if ${\bf a}\wedge{\bf b}$ is a positive multiple of ${\bf u}\wedge{\bf v}$.

\begin{theorem}\label{thm:simple}
If $|{\bf u}\wedge{\bf v}|=1$, then $\exp(\theta\,{\bf u}\wedge{\bf v})$ is simple counterclockwise rotation by $\theta$ radians in the oriented two--plane spanned by $({\bf u}$, ${\bf v})$.
\end{theorem}

\begin{proof}
By lemma \ref{lem:plane}, we may assume ${\bf u}$ and ${\bf v}$ are mutually orthonormal.  Using theorem \ref{thm:expsimple}, we then compute $\exp(\theta\,{\bf u}\wedge{\bf v}){\bf u}=\cos\theta\,{\bf u}-\sin\theta\,{\bf v}$, and $\exp(\theta\,{\bf u}\wedge{\bf v}){\bf v}=\sin\theta\,{\bf u}+\cos\theta\,{\bf v}$.  Moreover, if ${\bf w}$ is orthogonal to both ${\bf u}$ and ${\bf v}$, then $\exp(\theta\,{\bf u}\wedge{\bf v}){\bf w}={\bf 0}$.
\end{proof}

\subsection{Logarithm of a simple rotation}

From theorem \ref{thm:simple}, we see that every simple rotation $R$ can be expressed as the exponential of the wedge product of two vectors, say $R=\exp({\bf u}\wedge{\bf v})$.  In this case, we refer to ${\bf u}\wedge{\bf w}$ as a {\em logarithm} of $R$.  Logarithms of rotations are not unique, as the sine and cosine functions are periodic.  However, for simple rotations by angles less that $\pi$, there is a canonical choice for logarithm.

\begin{theorem}\label{thm:logsimple}
If $R$ is a simple rotation in dimension $n$, then the angle of rotation is given by $\theta=\cos^{-1}\tfrac{1}{2}\bigl({\it tr}(R)-n+2\bigr)$.  Moreover, if $0<\theta<\pi$, then $R=\exp(f)$, where $f=(\theta/2\sin\theta)(R-R^t)$.
\end{theorem}

\begin{proof}
From lemma \ref{lem:elemcomp}, ${\it tr}({\it proj}_{\bf uv})=2$ and ${\it tr}({\bf u}\wedge{\bf v})=0$.  Thus from theorem \ref{thm:expsimple}, ${\it tr}\bigl(\exp(\theta\,{\bf u}\wedge{\bf v})\bigr)=n-2(1-\cos\theta)$.
\end{proof}

If $0<\theta<\pi$, we write $\log(R)=f$.  If $\theta=0$, we simply take $\log(R)=I$.  And if $\theta=\pi$, then we are only able to deduce the projection onto the two--plane: ${\it proj}_{\bf uv}=\tfrac{1}{2}(I-R)$; $R$ is multiplication by $-1$ in this two--plane.

\section{Decomposing an antisymmetric matrix}
\label{sec:decomp}

In two and three dimensions, every rotation is simple.  In higher dimensions, one must compose simple rotations to obtain a general rotation.  If the two--planes of the constituent simple rotations are mutually orthogonal (from Lie Theory, this can always be done), we may obtain an explicit formula for the composite rotation.  We will do this in the next section, although we outline the general procedure here.

\begin{lemma}\label{lem:wedgeorth}
The two--planes spanned by the pairs ${\bf u}$, ${\bf v}$ and ${\bf a}$, ${\bf b}$ are orthogonal if and only if $({\bf u}\wedge{\bf v})({\bf a}\wedge{\bf b})=0$.
\end{lemma}

\begin{proof}
Since ${\bf u}\cdot{\bf a}={\bf u}\cdot{\bf b}={\bf v}\cdot{\bf a}={\bf v}\cdot{\bf b}=0$, lemma \ref{lem:elemcomp}\eqref{it:outercomp} implies that $({\bf u}\wedge{\bf v})({\bf a}\wedge{\bf b})=0$.  Conversely, if $({\bf u}\wedge{\bf v})({\bf a}\wedge{\bf b})=0$, then item lemma \ref{lem:elemcomp}\eqref{it:proj} implies that ${\it proj}_{\bf uv}{\it proj}_{\bf ab}=0$, and the two--planes are orthogonal.
\end{proof}

Suppose that $f_j$, $1\leq j\leq m$, are wedge products with mutually orthgonal two--planes, lemma \ref{lem:wedgeorth} then implies that $(\sum_jf_j)^k=\sum_jf_j^k$ for all $k>0$.  By summing the Taylor series, we thus find $\exp(\sum_jf_j)=I+\sum_j(\exp(f_j)-I)$.  Theorem \ref{thm:expsimple} can be then used to obtain an explicit formula for the latter.

If we are able to write an antisymmetric matix $f$ as a the sum of wedge products $f_j$ whose two--planes are mutually orthogonal, we will call the decomposition $f=\sum_jf_j$ an {\bf orthogonal decomposition} of $f$.  By dimensional considerations, there can be at most $\lfloor n/2\rfloor$ (where $\lfloor x\rfloor$ denotes the integer part of $x$) nozero summands in the orthogonal decomposition of a $n\times n$ antisymmetric matrix.

\subsection{Existence of orthogonal decomposition}

Up to this point, we have used only elementary considerations; however the the following lemma makes use of the {\em spectral theorem for real $n\times n$ antisymmetric matrices} $A$: there exists an orthogonal matrix $P$ such that $A=PBP^t$, where $B$ has $2\times 2$ blocks of the form $(\begin{smallmatrix}0 & \lambda\\-\lambda & 0\end{smallmatrix})$ along the diagonal; if $n$ is odd, $B$ has at least one additional row and column of zeros.  See \cite{Gantmacher}.

\begin{theorem}\label{thm:decomp0}
Every real antisymmetric matrix, is the sum of wedge products whose two--planes are mutually orthogonal.
\end{theorem}

\begin{proof}
Let $B$ and $Q$ be as above; and let ${\bf e}_j$ denote the $n$--dimensional vector whose $j$--th component is unity, and all other components are zero.  We may write $B=\sum_j\lambda_j{\bf e}_{j-1}\wedge{\bf e}_{j+1}$.  Thus $B$ is evidently the sum of wedge products with mutually orthogonal two--planes; and since $Q$ is orthogonal, $A=QBQ^T=\sum_j\lambda_j(Q{\bf e}_{j-1})\wedge(Q{\bf e}_{j+1})$ is also.
\end{proof}

\begin{corollary}\label{cor:decomp0}
Every real antisymmetric $n\times n$ matrix, with $n=4,5$, is the sum of at most two wedge products whose two--planes are orthogonal.\qed
\end{corollary}

\subsection{Finding an orthogonal decomposition}

\begin{theorem}\label{thm:decomp}
Let $f$ be a real $n\times n$ antisymmetric matrix, with $n=4,5$, then
\begin{equation*}
  \Delta\doteq{\it tr}(f^4)-\tfrac{1}{4}{\it tr}^2(f^2)
  \quad\text{and}\quad
  \theta_\pm^2\doteq-\tfrac{1}{4}{\it tr}(f^2)
                     \pm\tfrac{1}{2}\sqrt{\Delta}
\end{equation*}
are nonnegative.  Moreover if $\Delta\neq 0$, then we have the orthogonal decomposition $f=f_++f_-$ where
\begin{equation*}
  f_\pm=\mp\frac{1}{\sqrt\Delta}(\theta_\mp^2 f+f^3)
\end{equation*}
are wedge products with $f_+f_-=0$ and $|f_\pm|=\theta_\pm\doteq\sqrt{\theta_\pm^2}$.  This decomposition is unique up to reordering of summands.
\end{theorem}

\begin{proof}
By corollary \ref{cor:decomp0}, we may write $f$ as the sum of wedge products $f=f_++f_-$, with $f_+f_-=0$.  However from lemma \ref{lem:elemcomp}, $f_\pm^2=-|f_\pm|^2{\it proj}_\pm$, where ${\it proj}_\pm$ denotes orthogonal projection onto the two--plane of the wedge product $f_\pm$.  Therefore, we have $f^2=-|f_+|^2{\it proj}_+-|f_-|^2{\it proj}_-$; and by taking the trace, we get the equation ($\ast$) $|f_+|^2+|f_-|^2=-\tfrac{1}{2}\,{\it tr}(f^2)$.  Similarly, as the square of a projection is itself, one computes the trace of $f^4$ to obtain the equation: ($\ast\ast$) $|f_+|^4+|f_-|^4=\tfrac{1}{2}\,{\it tr}(f^4)$.  Setting $x\doteq|f_+|^2$ and $y\doteq|f_-|^2$, these two equations give us a system of two algebraic equations in the variables $x$, $y$.  Using ($\ast$) to eliminate $y$ in ($\ast\ast$), we arive at the quadratic equation $x^2+\tfrac{1}{2}{\it tr}(f^2)x+\tfrac{1}{8}{\it tr}^2(f^2)-\tfrac{1}{4}{\it tr}(f^4)=0$ with discriminant $\Delta$, and whose solutions are $x=\theta_\pm^2$.  With possible reindexing, we may take $x=\theta_+^2$ and $y=\theta_-^2$.  Note that ($\ast$) and ($\ast\ast$) also imply that $\Delta=2(|f_+|^4+|f_-|^4)-(|f_+|^2+|f_-|^2)^2=(|f_+|^2-|f_-|^2)^2\geq 0$.  In addition, $\theta_\pm^2=|f_\pm|^2\geq 0$.

To obtain the formulae for $f_\pm$, we use lemma \ref{lem:elemcomp}\eqref{it:cube} to get $f^3=f_+^3+f_-^3=-\theta_+^2f_+-\theta_-^2f_-$.  This leads to the Vandermonde system
\begin{equation*}
  \begin{pmatrix} 1          & 1\\
                  \theta_+^2 & \theta_-^2\end{pmatrix}
  \begin{pmatrix} f_+\\
                  f_-\end{pmatrix}
  =
  \begin{pmatrix} f\\
                 -f^3\end{pmatrix}.
\end{equation*}
The determinant of the $2\times 2$ matrix on the left is $\theta_-^2-\theta_+^2=-\sqrt\Delta$.  Since we are assuming $\Delta>0$, we may invert the system to get the stated (necessarily unique) solutions for $f_+$, $f_-$.
\end{proof}

In the case when $\Delta=0$ (so that $\theta_+^2=\theta_-^2$), an orthogonal decomposition exists, but will not be unique.  This stems from the fact that an eigenspace of dimension greater than unity has no natural choice of basis.

\section{Rotations in four and five dimensions}

As we mentioned previously, every rotation is necessarily the exponential of an antisymmetric matrix; see \cite{Curtis}.  We use this fact, coupled with orthogonal decomposition, to write down a parametrization for all rotations in four and five dimensions, and to find logarithms for such rotations.

\subsection{Rotations}

\begin{theorem}\label{thm:rot45}
Let $f$, $\Delta$, $\theta_\pm$ be as in theorem \ref{thm:decomp}.  If $\Delta\neq 0$ and $\theta_\pm\neq 0$, then
\begin{equation*}
  \exp(f)=I+\frac{1}{\sqrt\Delta}(Af+Bf^2+Cf^3+Df^4),
\end{equation*}
where
\begin{align*}
A&\doteq\frac{\theta_+^2\sin{\theta_-}}{\theta_-}
      -\frac{\theta_-^2\sin{\theta_+}}{\theta_+},
&B&\doteq\frac{\theta_+^2(1-\cos{\theta_-})}{\theta_-^2}
        -\frac{\theta_-^2(1-\cos{\theta_+})}{\theta_+^2},\\
C&\doteq\frac{\sin{\theta_-}}{\theta_-}
        -\frac{\sin{\theta_+}}{\theta_+},
&D&\doteq\frac{1-\cos{\theta_-}}{\theta_-^2}
        -\frac{1-\cos{\theta_+}}{\theta_+^2}.
\end{align*}
Moreover if $\Delta=0$ and $\theta_+\neq 0$, then $\displaystyle\exp(f)=I+\frac{\sin{\theta_+}}{\theta_+}f+\frac{1-\cos{\theta_+}}{\theta_+^2}f^2$.
\end{theorem}

\begin{proof}
By theorem \ref{thm:decomp} we have $f_+f_-=0$, so that $(\theta_-^2f+f^3)(\theta_+^2f+f^3)=0$; from which it follows that $f^6=-\theta_+^2\theta_-^2f^2-(\theta_+^2+\theta_-^2)f^4$.  If $\Delta\neq 0$ and $\theta_+\neq 0$, then lemma \ref{lem:elemcomp}\eqref{it:proj} implies ${\it proj}_+=-(1/\theta_+^2\Delta)(\theta_+^2f+f^3)^2=(1/\theta_+^2\Delta)(\theta_-^2f^2+f^4)$, and similarly for ${\it proj}_-$.  The first stated formula for $\exp(f)$ then follows from theorem \ref{thm:expsimple} and the discussion at the beginning of section \ref{sec:decomp}.  On the other hand, if $\Delta=0$ and $\theta_+\neq 0$, then $\theta_+=\theta_-$; so that $|f_+|=|f_-|$.  Thus $f^3=f_+^3+f_-^3=-|f_+|^2(f_++f_-)$ by lemma \ref{lem:elemcomp}\eqref{it:cube}; whence $f^3=-\theta_+^2f$.  We may then sum the Taylor series for $\exp(f)$ to obtain the second stated formula.
\end{proof}

The case $\Delta\neq 0$ and $\theta_+=0$ cannot occur, since $-{\it tr}(f^2)=|f|^2\geq 0$; if $\theta_-=0$, then $f$ is already a wedge product, and theorem \ref{thm:expsimple} applies.  There are no other nontrivial cases.  We remark that in the case when $n=4$, $\Delta\neq 0$, and $\theta_\pm\neq 0$, we can use the identity ${\it proj}_++{\it proj}_-=I$ to eliminate $f^4$ in the formula for $\exp(f)$.

\subsection{Logarithm of a rotation}

\begin{theorem}\label{thm:log45}
Let $R$ be a $n\times n$ rotation matrix, where $n=4,5$.  Define
\begin{equation*}
  \delta
  \doteq\tfrac{1}{2}{\it tr}(R^2)-\tfrac{1}{4}{\it tr}^2(R)
         +\tfrac{1}{2}(n-4){\it tr}(R)-\tfrac{1}{4}n(n-6)
\end{equation*}
\begin{equation*}
  y_\pm
  \doteq\tfrac{1}{4}\bigl({\it tr}(R)-n+4\bigr)
         \pm\tfrac{1}{2}\sqrt\delta.
\end{equation*}
Then $\delta\geq 0$ and $-1\leq y_\pm\leq 1$.  If $\delta\neq 0$ and $y_\pm\neq\pm 1$, then $\exp(f)=R$, where $f$ is antisymmetric and has the orthogonal decomposition $f=f_++f_-$ with
\begin{equation*}
  f_\pm\doteq\mp\frac{\theta_\pm}{2\sin\theta_\pm\sqrt\delta}
             \bigl(y_\mp(R-R^t)-\tfrac{1}{2}(R^2-R^{2t})\bigr),
\end{equation*}
where $\theta_\pm\doteq\cos^{-1}y_\pm$.  If $\delta=0$ and $y_+\neq\pm 1$, then $\displaystyle f=\frac{\theta_+}{2\sin\theta_+}(R-R^t)$.
\end{theorem}

\begin{proof}
We know that $R=\exp(\theta_+w_++\theta_-w_-)$, where $w_\pm$ are wedge products of unit length and whose two--planes are orthogonal, and where $0\leq\theta_\pm\leq\pi$.  By theorem \ref{thm:expsimple} and the discussion at the beginning of section \ref{sec:decomp},
\begin{equation}\label{eq:R}
  R=I-(1-\cos\theta_+){\it proj}_++\sin\theta_+w_+
     -(1-\cos\theta_-){\it proj}_-+\sin\theta_-w_-,
\end{equation}
where ${\it proj}_\pm=-w_\pm^2$.  Using lemma \ref{lem:elemcomp}\eqref{it:proj} and \eqref{it:cube}, orthogonality, and the double angle formulas $\sin{2\theta}=2\sin\theta\cos\theta$ and $\cos{2\theta}=\cos^2\theta-\sin^2\theta$, one computes
\begin{equation}\label{eq:R2}
\begin{split}
  R^2&=I-(1-\cos{2\theta_+}){\it proj}_++\sin{2\theta_+}w_+\\
     &\quad\quad\quad\quad
      -(1-\cos{2\theta_-}){\it proj}_-+\sin{2\theta_-}w_-.
\end{split}
\end{equation}
Set $y\doteq\cos\theta_+$ and $z\doteq\cos\theta_-$.  Since $\cos{2\theta}=2\cos^2\theta-1$, we get the equations ($\ast$) ${\it tr}(R)=2y+2z+n-4$ and ($\ast\ast$) ${\it tr}(R^2)=4y^2+4z^2+n-8$; and from these we get $y^2-\tfrac{1}{2}[{\it tr}(R)-n+4]y+\tfrac{1}{8}[{\it tr}^2(R)-{\it tr}(R^2)-2(n-4){\it tr}(R)+n^2-7n+8]=0$, whose discriminent is $\delta$, and solutions are $y=y_\pm$.  Moreover, one computes using ($*$) and ($**$) that $\delta=(y_+-y_-)^2$ and we may take $y_\pm=\cos\theta_\pm$.  Thus $\delta\geq 0$ and $-1\leq y_\pm\leq 1$.

Let $f_+\doteq\theta_+w_+$ and $f_-\doteq\theta_-w_-$.  First assume that $\delta\neq 0$ and $0<\theta_\pm<\pi$ (that is, $y_\pm\neq\pm 1$).  Equation \eqref{eq:R} then implies that $\tfrac{1}{2}(R-R^t)=\sin\theta_+f_+/\theta_++\sin\theta_-f_-/\theta_-$, and equation \eqref{eq:R2} implies that $\tfrac{1}{2}(R^2-R^{2t})=\sin{2\theta_+}f_+/\theta_++\sin{2\theta_-}f_-/\theta_-$.  Using the double angle formula for sine, we may recast these equations as the matrix equation
\begin{equation*}
  \begin{pmatrix}
    1    & 1\\
    2y_+ & 2y_-
  \end{pmatrix}
  \begin{pmatrix}
    \sin\theta_+f_+/\theta_+\\
    \sin\theta_-f_-/\theta_-
  \end{pmatrix}
  =
  \begin{pmatrix}
    (R-R^t)/2\\
    (R^2-R^{2t})/2
  \end{pmatrix}.
\end{equation*}
The determinant of the $2\times 2$ matrix on the left is $-2\sqrt\delta\neq 0$, so we can invert to get the first solution stated.  Second, if $\delta=0$, then $y_+=y_-$.  Therefore, equation \eqref{eq:R} becomes $R=I-(1-y_+)({\it proj}_++{\it proj}_-)+\sin\theta_+f/\theta_+$, which implies the second solution, provided that $0<\theta_+<\pi$.
\end{proof}

In the case when $\delta\neq 0$ and $y_+=1$ or $y_-=1$, $R$ is a simple rotation, and is handled by theorem \ref{thm:logsimple}.  Note that the case $\delta\neq 0$ and $y_+=-1$ is not possible (otherwise $y_-<-1$).  Moreover, from equation \eqref{eq:R}, we see that if $y_\pm=-1$ (so that $\delta=0$), then $R$ is multiplication by $-1$ in the four--plane with orthogonal projection map $(I-R)/2$ --- although this is only meaningful in the case $n=5$.  The only other nontrivial case occurs when $\delta\neq 0$ and $y_-=-1$; in which case $R=\exp(f_+)-2{\it proj}_-$, where $f_+=(\theta_+/2\sin\theta_+)(R-R^t)$ and ${\it proj}_-=(\exp(f_+)-R)/2$.

\section{Concluding remarks}

The formulas given for rotations and their logarithms in three dimensions seem to be well--known, although surprisingly obscure --- most likely due to the predominant use of Euler angles.  On the other hand, I have not seen the formulas given here for four and five dimensions previously in print, though it strikes me as unlikey that they are original.  Quaternionic formulas for the two--planes and rotation angles for a four--dimensional rotation, however, can be found in \cite{Weiner}.

The constructions given for the orthogonal decomposition, exponential map, and logarithm map may be extended to higher dimensions.  However for the two--plane rotation angles in dimension $n$, one must solve a polynomial of degree $\lfloor n/2\rfloor$, which cannot be done analytically for $n\geq 10$.  Moreover to find the orthogonal decomposition, one needs to invert a $\lfloor n/2\rfloor\times\lfloor n/2\rfloor$ matrix.


\end{document}